\documentclass[11pt]{article}
\newcommand{\proof}{\noindent {\bf Proof: }}

\newtheorem{lemma}{Lemma}
\newtheorem{statement}{Statement}

\def\qed{\hfill $\Box$}

\usepackage{amssymb}
\usepackage{graphicx}

\begin{document}
\title{Triangle in a brick}
\author{\'A.G.Horv\'ath\\ Department of Geometry, \\
Budapest University of Technology and Economics,\\
H-1521 Budapest,\\
Hungary}
\date{September 15, 2010}

\maketitle

\begin{abstract}
In this paper we shall investigate the following problem: Which is the largest regular triangle in a brick with volume 1 and having edge lengthes greater than $\frac{1}{\sqrt{2}}$? We prove that there are three optimal cases in all which the edge lengthes of the triangles are $\sqrt{2}$, respectively.
\end{abstract}

{\bf MSC(2000):} 52C17

{\bf Keywords:} brick with unit volume, regular triangle with maximal edge lenght

\section{Introduction and the results}

In the paper \cite{gho 1} we can find a lemma without complete (and exact) proof (see Lemma in \cite{gho 1}) was very important in the solving of the problem; find the optimal ball packing with respect to the crystallographic group $P2_12_12_1$. (We use the standard notation of a crystallographic group, it can be found e.g. in the International Tables \cite{inttable}.) In this paper we give a short and elegant proof for this basic statement, namely we shall prove that:
\begin{statement}
If a brick has a unit volume, its side lengths are not less than $\frac{1}{\sqrt{2}}$ and the vertices of a triangle lie on the skew edges of this brick, then there is a side of the triangle whose length is not greater than $\sqrt{2}$.
\end{statement}
This geometric problem can be handing by the standard analytic methods badly, we have to investigate certain complete polynomial equations of order six or of order higher than six. So without further geometric simplification there is no chance to solve it exactly.  The geometric intuition helps us to reduce the number and order of the variables thus it is very important in the solution of such a problem. Similar elementary but hard results can be found in the papers \cite{gho 2} and \cite{gho 3} also showing that the intuition cannot be change to a mechanical (theoretical or analytical) computation in every time and every cases.

The following simple lemma plays an important role in the proof of the required result.

\begin{figure}[ht]
  \centering
  \includegraphics[scale=0.5]{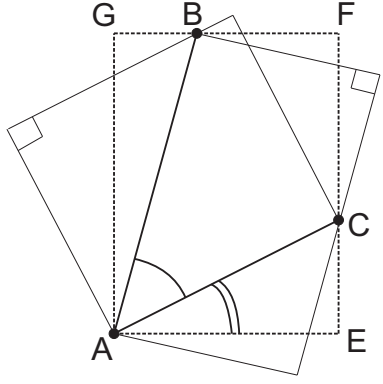}\\
  \caption{Extremal rectangles in the lemma.}
\end{figure}

\begin{lemma}
Let $A,B$ and $C$ three fixed point of the plane does not lying on a line and assume that the angle $BAC\angle$ is not obtuse. Then from among the rectangles which one vertex is $A$ and an edge of it opposite to $A$ and has minimal area one of its edges either $AB$ or $AC$.
\end{lemma}

On Fig. 1 we can see the possible rectangles with possible minimal areas.
From this lemma first we shall prove the following statement:

\begin{statement}
If a brick has a unit volume, its side lengths are not less than $\frac{1}{\sqrt{2}}$ and the vertices of a regular triangle lie on the skew edges of this brick, then the length of its sides is not greater than $\sqrt{2}$.
\end{statement}

Finally we observe that Statement 1 already is a simple consequence of Statement 2.
We note that the condition on the lengthes of the edges of the brick may not omit, in a sufficiently thin brick with volume 1 we can put whatever large regular triangle.

\section{Proofs}

\proof[Lemma 1] First we note that booth of the mentioned rectangles are exists because of the initial assumptions. The notation of Figure 1 we now have
$$
|AE|=|AC|\cos CAE\angle
$$
and
$$
|AG|=|AB|\sin (BAC\angle +CAE\angle)
$$
thus the examined area is
$$
T=|AB| |AC|\cos CAE\angle\sin (BAC\angle +CAE\angle).
$$
$T$ is minimal if
$$
f( CAE\angle)=\cos CAE\angle\sin (BAC\angle +CAE\angle)
$$
is minimal on the parameter domain defined by the inequalities $0<BAC\angle\leq \frac{\pi}{2}$ and $0\leq CAE\angle\leq \frac{\pi}{2}-(BAC\angle )$. The computation
$$
f( CAE\angle)=\cos CAE\angle\sin (BAC\angle +CAE\angle)=\frac{1}{2}\left[\sin (BAC\angle +2CAE\angle)+\sin (BAC\angle )\right]\geq
$$
$$
\geq \sin (BAC\angle )
$$
shows that it would be minimal in cases $CAE\angle=0$ or $CAE\angle=\frac{\pi}{2}-(BAC\angle )$
This proves the lemma.
\qed

Now we can prove the statement on the regular triangle.

\proof[Statement 2] We distinguish two cases:
\begin{enumerate}
\item There is an edge of the triangle (say $AB$) which is also an edge of the brick,
\item there is no such an edge.
\end{enumerate}

{\bf In the first case} the third vertex $C$ is the midpoint of the opposite edge of the brick. Let now $a,b,c,d$ the edge lengthes of the brick and the triangle, respectively. (See on Fig. 2.)

\begin{figure}[ht]
  \centering
  \includegraphics[scale=0.75]{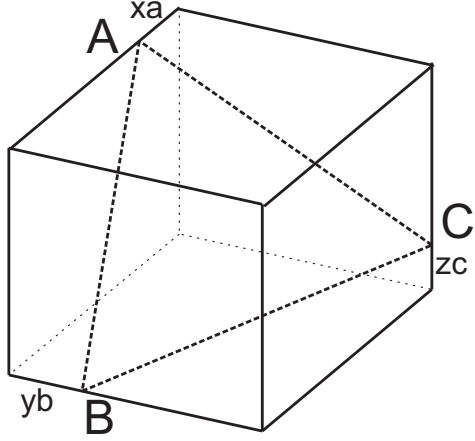}\\
  \caption{The case of regular triangle.}
\end{figure}

This implies the following equalities:
$$
d^2=a^2+b^2+\frac{1}{4}c^2=c^2
$$
and
$$
a^2+b^2=\frac{3}{4}c^2=\frac{3}{4}\frac{1}{a^2b^2}.
$$
Thus a short computation shows that
$$
d^2=\frac{2}{3}\left[b^2+\sqrt{b^4+\frac{3}{b^2}}\right].
$$
With respect to the various $t=b^2$ we can see immediately that of the function
$$
f(t)=\frac{2}{3}\left[t+\sqrt{t^2+\frac{3}{t}}\right]
$$
there is precisely one minima on the half line $[\frac{1}{2},\infty]$ in an inner point of the interval $[\frac{1}{2},1]$. On the other hand if $t=b^2>1$ (and by assumption $a^2\geq \frac{1}{2}$) then
$$
c^2>\frac{4}{3}(1+\frac{1}{2})=2
$$
implying that $a^2<\frac{1}{2}$ which is a contradiction. So we can assume that $t\in [\frac{1}{2},1]$ and thus we have to take the maximum of $f$ on this closed interval. They are attained on the boundary points of $[\frac{1}{2},1]$. Hence we get that in this case $d^2\leq 2$ and the optimum given by two configuration which are equivalent to each other. The parameter values:
$$
a^2=\frac{1}{2}, b=1, c^2=2
$$
and
$$
a=1, b^2=\frac{1}{2}, c^2=2.
$$

We now investigate {\bf the second case}, we assume that there is no edge of the triangle which is also an edge of the brick. First we note that an triangle which can pack (satisfying the conditions of the arrangement) in a brick with less volume then 1 cannot be maximal one. Consider that edge of the brick which contains the vertex $A$ (see on Fig. 2). Project orthogonally the arrangement onto one of the face of the brick orthogonal to this edge. Observe that the image of the regular triangle is a triangle with vertices $A',B',C'$ which is drawing in a rectangle on the manner described in Lemma 1. If the edges $A'B'$ and $A'C'$ are not edges of the rectangle, then (by Lemma 1) there is a rectangle with smaller area and the brick building up this base with the original height also contains the original regular triangle on its skew edges. This means that the triangle is not maximal. This can be saying too in the case of the vertices $B$ and $C$ and thus we get that one of the edges of the triangle is an diagonal of a face of the brick. (See Fig. 3)

\begin{figure}[ht]
  \centering
  \includegraphics[scale=0.75]{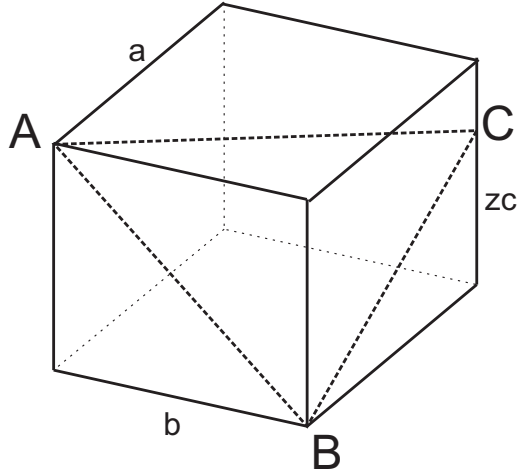}\\
  \caption{One edge is a diagonal.}
\end{figure}

By the notation of Fig.2 our problem analytically now simplified the following one: Minimalize the function:
$$
d^2=b^2+c^2
$$
where
$$
b^2+c^2=a^2+z^2c^2=a^2+b^2+(1-z)^2c^2
$$
and thus $z=\frac{b^2+c^2}{2c^2}$ and hence
$$
b^2+c^2=\frac{1}{b^2c^2}+\frac{(b^2+c^2)^2}{4c^2}.
$$
From this we get that
$$
c^2=\frac{1}{3}\left[-b^2+\sqrt{4b^4+\frac{12}{b^2}}\right]
$$
and
$$
d^2=\frac{2}{3}\left[b^2+\sqrt{b^4+\frac{3}{b^2}}\right].
$$
We can follow now the argument of the first case as we also get that $d^2\leq 2$. We now have two new optimal arrangements with the parameters
$$
b^2=\frac{1}{2}, c^2=\frac{3}{2}, a^2=\frac{4}{3};
$$
and
$$
a=b=c=1,
$$
respectively. All of the optimal cases can be seen on Fig. 4 and the statement is proved. \qed
\begin{figure}[ht]
  \centering
  \includegraphics[scale=1]{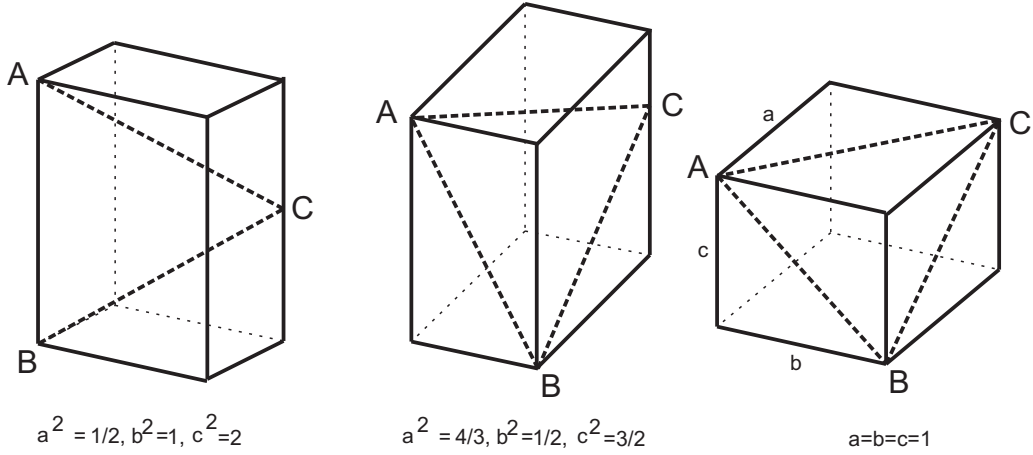}\\
  \caption{The maximal arrangements.}
\end{figure}

Finally we prove Statement 1.

\proof[Statement 1] The vertices of the triangle let denoted by $A,B,C$ and assume that $|AB|\leq |BC|\leq |AC|$. We distinguish two cases again.

\begin{enumerate}
\item The first inequality is a strictly one.
\item $|AB|=|BC|$
\end{enumerate}

{\bf In the first case} we have two possibilities because $A$ and $B$ are respective vertices of the cube. In fact, if $A$ is not a vertex then the strict inequalities $|AB|<|BC|,|AC|$ implies that the minimal edge of the triangle can be increased and it is not optimal triangle. So we have two cases
\begin{itemize}
                   \item The edge $AB$ is an edge of the brick. The value $|BC|$ is maximal if $C$ is the midpoint of the opposite edge. But if
                       $2<|AB|^2=a^2$ then $b^2c^2<\frac{1}{2}$ showing that $b^2$ or $c^2$ less than $\frac{1}{2}$. This is a contradiction.
                   \item The edge $AB$ is a diagonal of the brick. At this time $|BC|$ has maximal value if $C$ is the point of intersection of the bisector of the segment $AB$ and an edge of the face opposite to the face containing the diagonal $AB$. Now we have $$
                       \frac{1}{2}\leq |AB|<|BC|=|AC|
                       $$
                       contradicting Statement 2. In fact the affinity orthogonal to the face of $AB$ with ratio less than 1 gives a regular triangle of edge length $\sqrt{2}$ in a brick of volume less than 1. This is a contradiction again.
\end{itemize}
So in this case there is no optimal arrangement.

{\bf In the second case} we assume that $|AB|=|BC|\leq |AC|$. If $|BC|<|AC|$ and one of the vertices $A$ and $C$ is not a vertex of the brick then a little move of it can imply one of the inequalities $|AB|<|BC|\leq |AC|$ or $|BC|<|AB|\leq |AC|$ and thus the condition of the first case is holding. Thus we can assume that $A$ and $C$ are vertices of the brick and can apply the argument of the first case. So the only possibility to get an optimal arrangement if we assume that $|AB|=|BC|=|AC|$ and we get back the statement of Statement 2.
\qed


\begin{thebibliography}{99}

\bibitem{gho 1} \'A.G.Horv\'ath and E.Moln\'ar, Densest ball packing by orbits of the 10 fixed point free Euclidean space groups. (common with E.Moln\'ar) {\em Studia Sci. Math. Hungarica. } {\bf 29} (1994), 9--23.

\bibitem{gho 2} \'A.G.Horv\'ath and I.Prok, Packing congruent bricks in a cube. {\em J.
for Geometry and Graphics} {\bf 5} (2001), no.1, 1-11.

\bibitem{gho 3}  \'A.G.Horv\'ath, Maximal convex hull of connecting simplices. {\em Studies of the University of Zilina} {\bf 22} (2008), 7--19

\bibitem{inttable}
{\em International tables for X-ray crystallography.} Vol. A, Ed. by Theo Hahn, Reidel, 1983.



\end{thebibliography}
\end{document}